\documentclass[12pt,twoside, reqno]{article} \textheight=8.5truein
\usepackage{amsmath}
\usepackage{amssymb}
\usepackage{graphicx}

\usepackage{mathrsfs}
\usepackage{pb-diagram}

\def\stretchx{\Bumpeq{\!\!\!\!\!\!\!\!{\longrightarrow}}}

\newtheorem{theorem}{Theorem}[section]

\newtheorem{definition}{Definition}[section]
\newtheorem{remark}{Remark}[section]
\numberwithin{equation}{section}



\begin{document}

\title{{Chaotic dynamics in three dimensions:\\ a topological proof for a triopoly game model} 
}


\author{\textbf{Marina Pireddu} \\
Dipartimento di Matematica e Applicazioni, \\
Universit\`{a} degli Studi di Milano-Bicocca \\
via Cozzi 53, 20125 Milano\\
e-mail: marina.pireddu@unimib.it}

\maketitle

\begin{abstract}
\noindent 
We rigorously prove the existence of chaotic dynamics for the triopoly game model already studied, mainly from a numerical viewpoint, 
in \cite{NaTr->}. In the model considered, the three firms are heterogeneous and in fact each of them adopts a different decisional mechanism, i.e., linear approximation, best response and gradient mechanisms, respectively.\\
The method we employ is the so-called ``Stretching Along the Paths'' (SAP) technique in \cite{PiZa-07}, based on the Poincar\'e-Miranda Theorem and on the properties of the cutting surfaces.

\end{abstract}

\noindent \textbf{Keywords:} Chaotic dynamics; Stretching along the paths; triopoly games; heterogeneous players.

\vspace{2mm}

\noindent \textbf{2000 AMS subject classification:} 54H20, 54H25, 37B10, 37N40, 91B55.

\section{Introduction}\label{sec-a0}
$~$

In the economic literature, due to the complexity of the 
models considered, an analytical study of the associated dynamical features
turns out often to be too difficult or simply impossible
to perform. That is why many dynamical systems are studied mainly from a numerical viewpoint (see, for instance, \cite{Ag-99,BiTr-10,Tretal-10,YoMaPo-00}). 
Sometimes, however, even such kind of study turns out to be problematic, especially with high
dimensional systems, where several variables are involved.\\
In particular, as observed in Naimzada and Tramontana's working paper \cite{NaTr->}, this may be the reason for
the relatively low number of works on triopoly games (see, for instance, \cite{ElAgEl-09,Pu-96,TrEl-12}), where the context is given by an oligopoly composed
by three firms. In such framework, a local analysis can generally be performed in the
special case of homogeneous triopoly models, i.e., those in which the equations describing the dynamics are
symmetric (see, for instance, \cite{Ag-98, AgGaPu-00, RiSt-04}). \\
A more difficult task is that of studying heterogeneous triopolies, where instead the three
firms considered behave according to different strategies. 
This has been done, for instance, in \cite{Eletal-07, Ji-09}, as well as in the above mentioned paper by Naimzada and Tramontana \cite{NaTr->} where,
in addition to the classical heterogeneity with interacting agents adopting gradient and best response mechanisms, it is assumed that one of the firms adopts a linear approximation mechanism, which means that the firm does not know
the shape of the demand function and thus builds a conjectured demand function through the local knowledge of the true demand function. 
In regard to such model, those authors perform a stability analysis of the Nash equilibrium and show numerically that,
according to the choice of the parameter values, it undergoes a flip bifurcation or a Neimark-Sacker bifurcation leading to chaos.\\
What we then aim to do in the present paper is complementing that analysis, by proving the existence of 
chaotic sets only via topological arguments. This task will be performed using the ``Stretching Along the Paths'' (from now on, SAP) technique, already employed in \cite{MePiZa-09} to rigorously prove the presence of chaos for some discrete-time one- and bidimensional economic models of the classes of overlapping generations and duopoly game models. Notice however that, to the best of our knowledge, this is the first three-dimensional discrete-time application of the SAP technique, called in this way because it concerns maps that expand the arcs along one direction. 
We stress that, differently from other methods for the search of fixed points and the detection of chaotic dynamics based on more sophisticated algebraic or geometric tools, such as the Conley index or the Lefschetz number (see, for instance, \cite{Ea-75,MiMr-95,SrWo-97}), the SAP method relies on relatively elementary arguments and it is easy to apply in practical contexts, without the need of ad-hoc constructions. No differentiability conditions are required for the map describing the dynamical system under analysis and even continuity is needed only on particular subsets of its domain. Moreover, the SAP technique can be used to rigorously prove the presence of chaos also for continuous-time dynamical systems. In fact, in such framework it suffices to apply the results in Section \ref{sec-a1}, suitably modified, to the Poincar\'e map associated to the considered system and thus one is led back to work with a discrete-time dynamical system. However, the geometry required to apply the SAP method turns out to be quite different in the two contexts: in the case of discrete-time dynamical systems we look for ``topological horseshoes'' (see, for instance, \cite{BuWe-95, KeYo-01, ZgGi-04}), that is, a weaker version of the celebrated Smale horseshoe in \cite{Sm-65}, while in the case of continuous-time dynamical systems one has to consider the case of switching systems and the needed geometry is usually that of the so-called ``Linked Twist Maps'' (LTMs) (see \cite{BuEa-80, De-78, Pr-86}), as shown for the planar case in \cite{PaPiZa-08,PiZa-08}. We also stress that the Poincar\'e map is a homeomorphism onto its
image, while in the discrete-time framework the function describing the considered dynamical system need not be one-to-one, like in our example in Section \ref{sec-a2}. Hence, in the latter context, it is in general not be possible to apply the results for the Smale horseshoe, where one deals with homeomorphisms or diffeomorphisms.
As regards three-dimensional continuous-time applications of the SAP method, those have recently been performed in \cite{RHZa-12}, in a higher-dimensional counterpart of the LTMs framework, and in \cite{RHZa->}, where a system switching between different regimes is considered.\\
For the reader's convenience, we are going to recall in Section \ref{sec-a1} what are the basic mathematical ingredients behind the SAP method, as well as the main conclusions it allows to draw about the chaotic features of the model under analysis.
It will then be shown in Section \ref{sec-a2} how it can be applied to the triopoly game model taken from \cite{NaTr->}.
Some further considerations and comments can be found in Section \ref{sec-a3}, which concludes the paper.

\section{The ``Stretching along the paths'' method}\label{sec-a1}
$~$

In this section we briefly recall what the ``Stretching along the paths'' (SAP) technique consists in,
referring the reader interested in further mathematical details to \cite{PiZa-07}, where the original planar theory
by Papini and Zanolin in \cite{PaZa-04a,PaZa-04b} has been extended to the $N-$dimensional setting, with $N\ge 2.$ \\
In the bidimensional setting, elementary
theorems from plane topology suffice, while in the higher-dimensional framework some
results from degree theory are needed, leading to the study of the so-called
``cutting surfaces''.
In fact, the proofs of the main results in \cite{PiZa-07} (and in particular of Theorem \ref{th-fp} below), we do not recall here, are based on the properties of the cutting surfaces and on the Poincar\'e-Miranda Theorem, that is, an $N$-dimensional version of the Intermediate Value Theorem. \\ 
Since in Section \ref{sec-a2} we will deal with the three-dimensional setting only, we directly present the theoretical results in the special case in which $N=3.$
 
We start with some basic definitions.\\
A \textit{path} in a metric space $X$ is a continuous map $\gamma:
[t_0,t_1]\to X.$ We also set $\overline{\gamma}:=\gamma([t_0,t_1]).$
Without loss of generality, we usually take the unit interval
$[0,1]$ as the domain of $\gamma.$ A \textit{sub-path} $\sigma$ of
$\gamma$ is the restriction of $\gamma$ to a compact sub-interval
of its domain. 
By a \textit{generalized parallelepiped} 
we mean a set ${\mathcal P}\subseteq X$ which is homeomorphic to the
unit cube $I^3:=[0,1]^3,$ through a homeomorphism $h: {\mathbb R}^3\supseteq I^3 \to \mathcal P\subseteq X.$
We also set
$${\mathcal P}^{-}_{\ell}:= h([x_3 = 0])\,,\quad
{\mathcal P}^{-}_{r}:= h([x_3 = 1])$$ and call them the \textit{left} and
the \textit{right} faces of $\mathcal P,$ respectively,
where\footnote{Notice that the choice of privileging the third coordinate is purely conventional. In fact, any other choice would give the same results, as it is possible to compose the homeomorphism $h$ with a suitable permutation on three elements, without modifying its image set.}
$$[x_3 = 0]:= \{(x_1,x_2,x_3)\in I^3: \, x_3 = 0\} \,\,\,
\mbox{ and } \,\,\, [x_3 = 1]:= \{(x_1,x_2,x_3)\in I^3: \, x_3 = 1\}.$$ 
Setting
$$\mathcal P^{-}:= \mathcal P^{-}_{\ell}\cup \mathcal P^{-}_{r}\,,$$
we call the pair
$${\widetilde{\mathcal P}}:= (\mathcal P, \mathcal P^-)$$
an {\textit{oriented parallelepiped}} \textit{of $X$}.

\medskip

Although in the application discussed in the present paper the
space $X$ is simply $\mathbb R^3$ and the
generalized parallelepipeds are standard parallelepipeds, the generality of our definitions 
makes them applicable in different contexts (see Figure 1).

\smallskip

We are now ready to introduce the {\it stretching along the paths} property for maps between oriented rectangles.

\begin{definition}[SAP]\label{def-sap}
{\rm{Let ${\widetilde{\mathcal A}}:=
({\mathcal A},{\mathcal A}^-)$ and ${\widetilde{\mathcal B}}:=
({\mathcal B},{\mathcal B}^-)$ be oriented parallelepipeds of a metric
space $X.$ Let also $\psi: \mathcal A\to X$ be a function
and ${\mathcal K}\subseteq {\mathcal A}$
be a compact set. We say that \textit{$({\mathcal K},\psi)$ stretches
${\widetilde{\mathcal A}}$ to ${\widetilde{\mathcal B}}$ along the
paths}, and write
\begin{equation}\label{eq-sap}
({\mathcal K},\psi): {\widetilde{\mathcal A}} \stretchx {\widetilde{\mathcal B}},
\end{equation}
if the following conditions hold:
\begin{itemize}
\item{} \; $\psi$ is continuous on ${\mathcal K}\,;$ 
\item{} \; for every path $\gamma: [0,1]\to {\mathcal A}$ with
$\gamma(0)$ and $\gamma(1)$ belonging to different components of ${\mathcal A}^-,$ there exists a sub-path 
$\sigma:=\gamma|_{[t',t'']}:[0,1]\supseteq [t',t'']\to {\mathcal K},$ such that
$\psi(\sigma(t))\in {\mathcal B},\,
\forall\, t\in [t',t''],$ and, moreover, $\psi(\sigma(t'))$ and
$\psi(\sigma(t''))$ belong to different components of ${\mathcal B}^-.$
\end{itemize}
}}
\end{definition}

For a description of the relationship between the SAP relation and other ``covering relations'' in the literature on expansive-contractive maps, we refer the interested reader to \cite{MePiZa-09}.\\
A first crucial feature of the SAP relation is that, when it is satisfied with ${\widetilde{\mathcal A}}={\widetilde{\mathcal B}}$ \footnote{Note that this means both that $\mathcal A$ and $\mathcal B$ coincide as subsets of $X$ and that they have the same orientation. In fact, it is easy to find counterexamples to Theorem \ref{th-fp} if the latter property is violated (see, for instance, \cite{Pi-09}, pag. 11).}, it ensures the existence of a fixed point localized in the compact set $\mathcal K.$ In fact the following result does hold true.

\begin{theorem}\label{th-fp}
Let ${\widetilde{\mathcal P}}:= ({\mathcal P},{\mathcal P}^-)$
be an oriented parallelepiped of a metric space $X$ and let $\psi: \mathcal P\to X$ be a function. 
If ${\mathcal K}\subseteq {\mathcal P}$
is a compact set such that
$$({\mathcal K},\psi): {\widetilde{\mathcal P}} \stretchx {\widetilde{\mathcal P}},$$
then there exists at least a point $z\in {\mathcal K}$ with $\psi(z) = z.$
\end{theorem}

For a proof, see \cite{PiZa-07}, pagg. 307-308. Notice that the arguments employed therein are different from the ones used to prove the same result in the planar context (see, for instance, \cite{MePiZa-09}, pagg. 3301-3302), which are in fact much more elementary.\\
A graphical illustration of Theorem \ref{th-fp} can be found in Figure 1, where it looks evident that, differently from the classical Rothe and Brouwer Theorems, we do not require that $\psi(\partial \mathcal A)\subseteq \mathcal A.$

\begin{figure}[ht]\label{fig-0}
\centering
\includegraphics[scale=0.45]{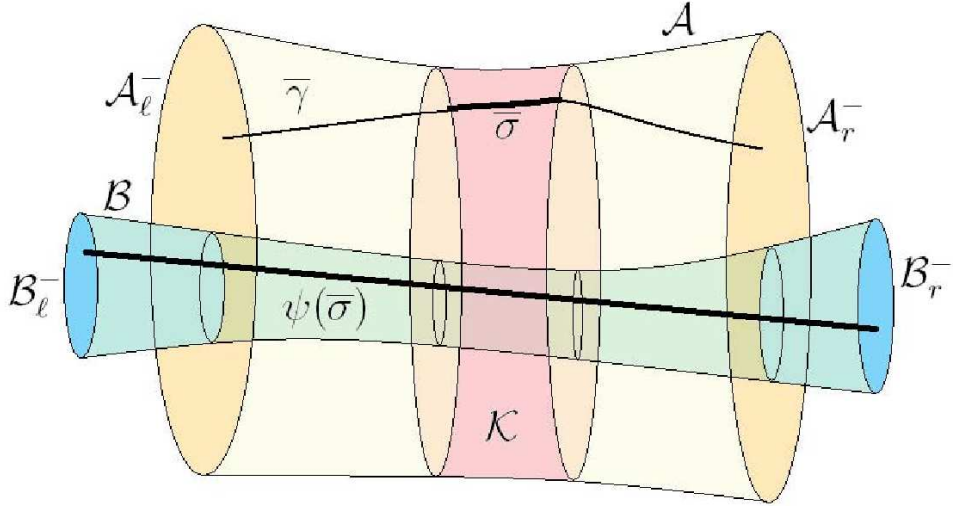}
\caption{The tubular sets $\mathcal A$ and $\mathcal B$ in the picture
are two generalized parallelepipeds, for which we
have put in evidence the compact set $\mathcal K$ and the boundary
sets $\mathcal A_{\ell}^-$ and $\mathcal A_{r}^-,$ as well as $\mathcal B_{\ell}^-$ and $\mathcal B_{r}^-.$ In this
case $(\mathcal K,\psi)$ stretches the paths of $\mathcal A$ across $\mathcal A$ itself and therefore the existence of a fixed
point for $\psi$ in $\mathcal K$ is ensured by Theorem
\ref{th-fp}.} 
\end{figure}

\smallskip

The most interesting case in view of detecting chaotic dynamics is when there exist pairwise disjoint compact sets playing the role of ${\mathcal K}$ in Definition \ref{def-sap}. Indeed, applying Theorem \ref{th-fp} with respect to each of them, we get a multiplicity of fixed points localized in those compact sets.
Another crucial property of the SAP relation is that it is preserved under composition of maps, and thus, when dealing with the iterates of the function under consideration, it allows to detect the presence of periodic points of any period (see Lemma A.1, Theorems A.1 and A.2 in \cite{MePiZa-09}, which can be directly transposed to the three-dimensional setting, with the same proof).\\
We now describe in Definition \ref{def-ch} what we mean when we talk about ``chaos'' and we explain in Theorem \ref{th-ch} which is the relationship between that concept and the stretching relation in Definition \ref{def-sap}.
We stress that Theorem \ref{th-ch} is the main theoretical result we are going to apply in Section \ref{sec-a2} and that it can be shown exploiting the two properties of the SAP relation mentioned above. In fact, its proof follows by the same arguments in Theorems 2.2 and 2.3 in \cite{MePiZa-09}.

\begin{definition}\label{def-ch}

\rm{Let $X$ be a metric space and let $\psi: X\supseteq D\to X$
be a function. Let also $m\ge 2$ be an integer and let
${\mathcal K}_0,\dots,{\mathcal K}_{m-1}$ be nonempty pairwise disjoint compact subsets of ${\mathcal D}.$ We say that
\textit{$\psi$ induces chaotic dynamics on $m$ symbols on the set
${\mathcal D}$ relatively to ${\mathcal K}_0,\dots,{\mathcal K}_{m-1}$} 
if, setting 
$$
\mathcal K:=\bigcup_{i=0}^{m-1}\mathcal K_i\subseteq {\mathcal D}
$$
and defining the nonempty compact set
\begin{equation}\label{eq-lam}
{\mathcal I}_{\infty}:=\bigcap_{n=0}^{\infty}\psi^{-n}(\mathcal K),
\end{equation}
then there exists a nonempty compact set
$${\mathcal I}\subseteq {\mathcal I}_{\infty} \subseteq {\mathcal K},$$
on  which the following conditions are fulfilled:

\begin{itemize}
\item[$(i)$] $\psi({\mathcal I}) = {\mathcal I};$

\item[$(ii)$] $\psi|_{\mathcal I}$ is semi--conjugate to the
Bernoulli shift on $m$ symbols, that is, there exists a continuous
map $\pi:{\mathcal I}\to\Sigma_m^+,$ where $\Sigma_m^+:=\{0,1\}^{\mathbb
N}$ is endowed with the distance
\begin{equation*}
\hat d(\textbf{s}', \textbf{s}'') := \sum_{i\in {\mathbb N}}
\frac{d(s'_i, s''_i)}{m^{i + 1}}\,,\quad \mbox{ for }\;
\textbf{s}'=(s'_i)_{i\in {\mathbb N}}\,,\,
\textbf{s}''=(s''_i)_{i\in {\mathbb N}}\in \Sigma_m^+
\end{equation*}
($\,d(\cdot\,,\cdot)$ is the discrete distance on $\{0,1\},$ i.e.,
$d(s'_i, s''_i)=0$ for $s'_i=s''_i$ and $d(s'_i, s''_i)=1$ for
$s'_i\not=s''_i$), such that the diagram
\begin{equation}\label{diag-1}
\begin{diagram}
\node{{\mathcal I}} \arrow{e,t}{\psi} \arrow{s,l}{\pi}
      \node{{\mathcal I}} \arrow{s,r}{\pi} \\
\node{\Sigma_m^+} \arrow{e,b}{\sigma}
   \node{\Sigma_m^+}
\end{diagram}
\end{equation}
commutes, where $\sigma:\Sigma_m^+\to\Sigma_m^+$ is the Bernoulli
shift defined by $\sigma((s_i)_i):=(s_{i+1})_i,\,\forall
i\in\mathbb N\,;$

\item[$(iii)$] the set of the periodic points of $\psi|_{{\mathcal I}_{\infty}}$ is dense in ${\mathcal I}$
and the pre--image $\pi^{-1}(\textbf{s})\subseteq {\mathcal I}$ of
every $k$-periodic sequence $\textbf{s} = (s_i)_{i\in {\mathbb N}}\in \Sigma_m^+$
contains at least one $k$-periodic point.
\end{itemize}
}
\end{definition}

\begin{remark}\label{cons}
According to Theorem 2.2 in \cite{MePiZa-09}, from $(ii)$ in Definition \ref{def-ch} it follows that:
\begin{itemize}
\item[$-$] $h_{\rm top}(\psi)\ge h_{\rm top}(\psi|_{\mathcal I})\geq h_{\rm top}(\sigma) = \log(m),$
where $h_{\rm top}$ is the topological entropy;

\item[$-$] there exists a compact invariant set
$\Lambda \subseteq {\mathcal I}$ such that $\psi\vert_{\Lambda}$
is semi--conjugate to the Bernoulli shift on $m$ symbols,
topologically transitive and displays sensitive dependence on initial
conditions.
\end{itemize}
\end{remark}

\begin{theorem}\label{th-ch}
Let ${\widetilde{\mathcal P}}:= ({\mathcal
P},{\mathcal P}^-)$ be an oriented parallelepiped of a metric space
$X$ and let $\psi: \mathcal P\to X$ be a function. If
${\mathcal K_0},\dots,{\mathcal K_{m-1}}$ are $m\ge 2$ pairwise disjoint compact
subsets of ${\mathcal P}$ such that
\begin{equation}\label{sap-rel}
({\mathcal K}_i,\psi): {\widetilde{\mathcal P}} \stretchx {\widetilde{\mathcal P}}, \mbox{ for } i=0,\dots,m-1,
\end{equation}
then $\psi$ induces chaotic dynamics on $m$ symbols on
${\mathcal P}$ relatively to ${\mathcal K}_0,\dots,{\mathcal K}_{m-1}.$
\end{theorem}

Notice that if the function $\psi$ in the above statement is also one--to--one on $\mathcal K:=\bigcup_{i=0}^{m-1}{\mathcal K}_i,$ then it is additionally possible to prove that $\psi$ restricted to a suitable invariant subset of $\mathcal K$ is semi--conjugate to the two--sided Bernoulli shift $\sigma:\Sigma_2\to\Sigma_2,$ $\sigma((s_i)_i):=(s_{i+1})_i,\,\forall i\in\mathbb Z,$ where $\Sigma_2:=\{0,1\}^{\mathbb Z}$ (see \cite[Lemma 3.2]{PiZa-08})\,\footnote{This is not the case in our application in Section \ref{sec-a2}. Indeed, as it looks clear from Figure 2, the map $F$ in \eqref{eq-tg} is not injective on the set $\mathcal K_0\cup\mathcal K_1$ introduced in Theorem \ref{th1}.}.
\medskip

We are now in position to explain what the SAP method consists in.
Given a dynamical system generated by a map $\psi,$ our technique consists in finding 
a subset $\mathcal P$ of the domain of $\psi$ homeomorphic to the unit cube and at least two disjoint compact subsets of $\mathcal P$ for which the stretching property in \eqref{sap-rel} is satisfied (when $\mathcal P$ is suitably oriented). In this way, Theorem \ref{th-ch} ensures the existence of chaotic dynamics in the sense of Definition \ref{def-ch} for the system under consideration and, in particular, the positivity of the topological entropy for $\psi,$ which is in fact generally considered as one of the trademark features of chaos.

\section{The triopoly game model}\label{sec-a2}
$~$

In this section we apply the SAP method to an
economic model belonging to the class of triopoly games, taken from \cite{NaTr->}.\\ 
By oligopoly, economists denote a market form characterized by
the presence of a small number of firms. Triopoly is a special
case of oligopoly where the firms are three. The term game
refers to the fact that the players - in our case the firms - make
their decisions reacting to each other actual or
expected moves, following a suitable strategy.
In particular, we will deal with a dynamic game where
moves are repeated in time, at discrete, uniform intervals.\\
More precisely, the model analyzed
can be described as follows.\\
The economy consists of three firms producing an identical commodity
at a constant unit cost, not necessarily
equal for the three firms. The commodity is sold in a single market
at a price which depends on total output through a given inverse
demand function, known to one firm (say, Firm $2$) globally and to another firm (say, Firm $1$) locally.
In fact, Firm $1$ linearly approximates the demand function around the latest realized pair of quantity and market price.
Finally, Firm $3$ does not know anything about the demand function and adopts a myopic adjustment mechanism, 
i.e., it increases or decreases its output according to the sign of the marginal profit from the last period.
The goal of each firm is
the maximization of profits, i.e., the difference between revenue
and costs. The problem of each firm is to decide at the beginning of every time period $t$
how much to produce in the same period on the basis of the limited
information available and, in particular, on the expectations about
its competitors' future decisions.\\
In what follows, we introduce the needed notation and the postulated assumptions:

\vskip .5cm

\noindent 1. \bf Notation   \rm \vskip .25cm

$x_t$: output of Firm 1 at time $t\,;$

$y_t$: output of Firm 2 at time $t\,;$

$z_t$: output of Firm 3 at time $t\,;$

$p$: unit price of the single commodity\,. \vskip .5cm

\noindent \bf  2. Inverse demand function 
\begin{equation}\label{p}
p:=\frac{1}{x+y+z}\,.
\end{equation}

\noindent \bf 3. Technology \rm \vskip .25cm

The unit cost of production for firm $i$ is equal to $c_i, \,i=1,2,3,$
where $c_1, c_2, c_3$ are (possibly different) positive constants. \vskip .5cm

\noindent \bf 4. Price approximation \rm \vskip .25cm

Firm 1 observes the
current market price $p_t$ and the corresponding total supplied quantity $Q_t=x_t+y_t+z_t.$ By
using market experiments, that player obtains the slope of
the demand function at the point $(Q_t, p_t)$ and, in the absence of other information, it
conjectures that the demand function, which has to pass through that point, is linear. 

\vskip .5cm

\noindent \bf 5. Expectations \rm \vskip .25cm

In the presence of incomplete information concerning their competitors'
future decisions (and therefore about future prices), Firms 1 and 2 are assumed to use naive expectations. This means that at each time $t$
both Firm 1 and 2 expect that the other two firms will keep output unchanged w.r.t. the 
previous period. \vskip .35cm

\medskip
\noindent
As shown in \cite{NaTr->}, the assumptions above lead to the following system of three difference equations in
the variables $x,\,y$ and $z$:

\begin{equation*}
\left\{
\begin{array}{ll}
x_{t+1}= \frac{2x_t+y_t+z_t-c_1(x_t+y_t+z_t)^2}{2}\\
\vspace{-3mm}\\
y_{t+1}= \sqrt{\frac{x_t+z_t}{c_2}}-x_t-z_t\\
\vspace{-3mm}\\
z_{t+1}= z_t+\alpha z_t\left(-c_3+\frac{x_t+y_t}{(x_t+y_t+z_t)^2}\right)\\
\end{array}
\right. \eqno{(\mbox{TG})}
\end{equation*}
where $\alpha$ is a positive parameter denoting the speed of Firm 3's adjustment to changes in profit and $c_1,\,c_2,\,c_3$ are the marginal costs.\\
We refer the interested reader to \cite{NaTr->} for a more detailed explanation of
the model, as well as for the derivation of (TG).

As mentioned in the Introduction, in \cite{NaTr->} Naimzada and Tramontana discuss the
equilibrium solution of system (TG) along with its stability and
provide numerical evidence of the presence of chaotic dynamics.
In particular, it is shown the existence of a double route 
to chaos: according to the parameter values, the Nash equilibrium
can undergo a flip bifurcation or a Neimark-Sacker bifurcation.
Moreover, in \cite{NaTr->} the authors numerically find multistability of different
coexisting attractors and identify their basins of attraction through a global analysis.\\
Hereinafter we will integrate that study
rigorously proving that, for certain parameter configurations,
system (TG) exhibits chaotic behavior in the precise sense
discussed in Section \ref{sec-a1} \footnote{Notice that, as we shall stress in Section \ref{sec-a3}, we
only prove \textit{existence} of an invariant, chaotic set, not
its \textit{attractiveness}.}.\\ 
In order to apply the SAP method
to analyze system (TG), it is expedient to represent it
in the form of a continuous map ${F=(F_1,F_2,F_3):\mathbb R_+^3\to\mathbb R^3},$ with components
\begin{equation}\label{eq-tg}
\begin{array}{ll}
{F_1(x,y,z):=\frac{2x+y+z-c_1(x+y+z)^2}{2}},\\
\vspace{-3mm}\\
{F_2(x,y,z):=\sqrt{\frac{x+z}{c_2}}-x-z},\\
\vspace{-3mm}\\
{F_3(x,y,z):=z+\alpha z\left(-c_3+\frac{x+y}{(x+y+z)^2}\right)}.\\
\end{array}
\end{equation}

We prove that the SAP property for the map $F$ is satisfied
when choosing a generalized rectangle in the family of parallelepipeds of the first quadrant described analytically by
\begin{equation}\label{eq-cu}
{\mathcal R=\mathcal R(x_i,\,y_i,\,z_i):=\left\{(x,y,z)\in\mathbb R^3: x_{\ell}\le x\le x_r,\,y_{\ell}\le y\le y_r,\,z_{\ell}\le z\le z_r\right\}},
\end{equation}
with $x_{\ell}<x_r,\,y_{\ell}<y_r,\,z_{\ell}<z_r$ and $x_i,\,y_i,\,z_i,\,i\in\{\ell,r\},$ satisfying the conditions in Theorem \ref{th1}.\\
The parallelepiped $\mathcal R$ can be oriented by setting
\begin{equation}\label{eq-or}
{\mathcal R^-_{\ell}:= [x_{\ell},x_r]\times [y_{\ell},y_r]\times\{z_{\ell}\}} \,  \mbox{ and } \, {\mathcal R^-_{r}:=[x_{\ell},x_r]\times [y_{\ell},y_r]\times\{z_r\}\,}.
\end{equation}

Consistently with \cite{NaTr->}, we choose the marginal costs as
$c_1=0.4,\, c_2=0.55$ and $c_3=0.6.$ On the other hand, in order to easily apply the SAP method we need the parameter $\alpha$ to be close to $17,$ while in \cite{NaTr->} the presence of chaos is numerically proven for $\alpha$ around $8$ \footnote{As explained below, it would be possible to apply our technique with a lower value for $\alpha,$ at the cost of changing the parameter conditions in Theorem \ref{th1} and of making the computations in the proof much more complicated. However, it seems not possible to apply the SAP method to the first iterate of $F$ when $\alpha$ is close to $8,$ which is the largest value considered in \cite{NaTr->}.}. The implications of this discrepancy will be discussed in Section \ref{sec-a3}.

\begin{figure}[ht]\label{fig-1}
\centering
\includegraphics[scale=0.7]{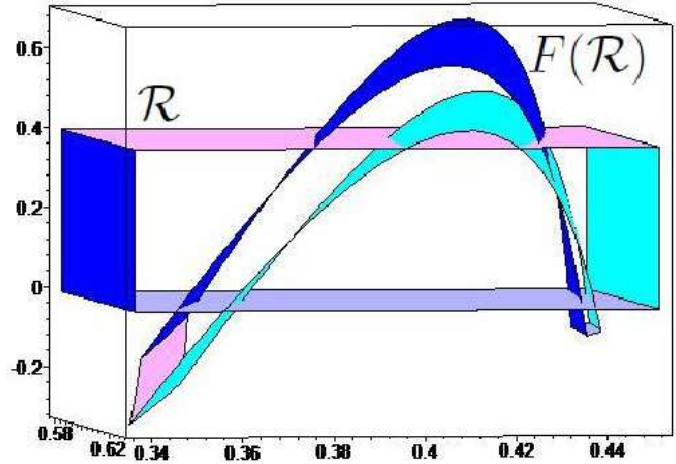}
\caption{ A possible choice of the parallelepiped $\mathcal
R$ for system (TG), according to conditions $(H1)$--$(H5).$ It has been oriented by
taking as $[\,\cdot\,]^{-}$-set the union of the two horizontal faces
$\mathcal R^-_{\ell}$ and $\mathcal R^-_{r}$ defined in
\eqref{eq-or}. In addition to $F(\mathcal R^-_{\ell})$ and $F(\mathcal R^-_{r}),$ we also
represent the image set of two vertical faces of $\mathcal R.$ Notice that we used the same color to depict a set and its $F$-image set.}
\end{figure}

\begin{figure}[ht]\label{fig-2}
\centering
\includegraphics[scale=0.7]{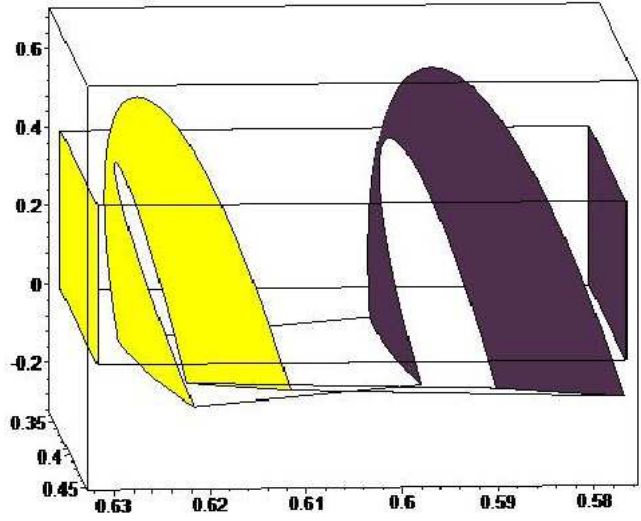}
\caption{This picture complements the previous one, by showing how the two vertical faces of $\mathcal R$ not considered in Figure 1 are transformed by the map $F.$ Again, the same color is used to depict a set and its $F$-image set.}
\end{figure}

\smallskip
\noindent Our result on system (TG) can be stated as follows:
\begin{theorem}\label{th1}
If the parameters of the map $F$ defined in \eqref{eq-tg} assume the following values 
\begin{equation}\label{par}
c_1=0.4,\,\, c_2=0.55,\,\,c_3=0.6,\,\,\alpha= 17,
\end{equation}
then, for any parallelepiped ${\mathcal R}=\mathcal R(x_i,y_i,z_i)$ belonging to the family described in \eqref{eq-cu}, with
$x_i,\,y_i,\,z_i,\,i\in\{\ell,r\},$ satisfying the conditions:
\begin{equation*}
\begin{array}{lll}
& {}\!\!\!\!\! (H1) & z_{\ell}=0\,;\\
& {}\!\!\!\!\! (H2) &  x_{\ell}+y_{\ell}>z_r\ge\sqrt{\frac{\alpha}{\alpha c_3-1}(x_{\ell}+y_{\ell})}-(x_{\ell}+y_{\ell})>0\,;\\
& {}\!\!\!\!\! (H3) & 2\left(\sqrt{\frac{\alpha}{\alpha c_3+1}(x_{r}+y_{r})}-(x_{r}+y_{r})\right)>z_r\,;\\
& {}\!\!\!\!\! (H4) &\frac{1}{c_1}-x_r>y_{r}+z_{r}>\frac{1}{2c_1}-x_{\ell}>0\,,\quad \frac{1}{2c_1}-x_r > y_{\ell}+z_{\ell}\,,\quad x_{r}\ge\frac{1}{4c_1}\,,\\
& & \frac{1}{2c_1}\left(1-c_1(y_{\ell}+y_{r}+z_{\ell}+z_{r})\right)\ge x_{\ell}>0\,,\quad \sqrt{\frac{y_{\ell}+z_{\ell}}{c_1}}-(y_{\ell}+z_{\ell})\ge x_{\ell}\,;\\
& {}\!\!\!\!\! (H5) &
x_{\ell}+z_{\ell}>\frac{1}{4c_2}\,,\quad y_r\ge\sqrt{\frac{x_{\ell}+z_{\ell}}{c_2}}-(x_{\ell}+z_{\ell})>0\,,\quad\sqrt{\frac{x_{r}+z_{r}}{c_2}}-(x_{r}+z_{r})\ge y_{\ell}>0\,,
\end{array}
\end{equation*}
and oriented as in \eqref{eq-or}, there exist two disjoint compact subsets ${\mathcal K}_0={\mathcal K}_0(\mathcal R)$ and ${\mathcal K}_1={\mathcal K}_1(\mathcal R)$ of $\mathcal R$ such that
\begin{equation}\label{eq-ks}
({\mathcal K}_i,F): {\widetilde{\mathcal R}} \stretchx {\widetilde{\mathcal R}}, \mbox{ for } i=0,1.
\end{equation}
Hence, the map $F$ induces chaotic dynamics on two symbols on
${\mathcal R}$ relatively to $\mathcal K_0$ and
$\mathcal K_1$ and displays all the properties listed in Theorem
\ref{th-ch}.
\end{theorem}

\bigskip

\begin{figure}[ht]\label{fig-3}
\centering
\includegraphics[scale=0.7]{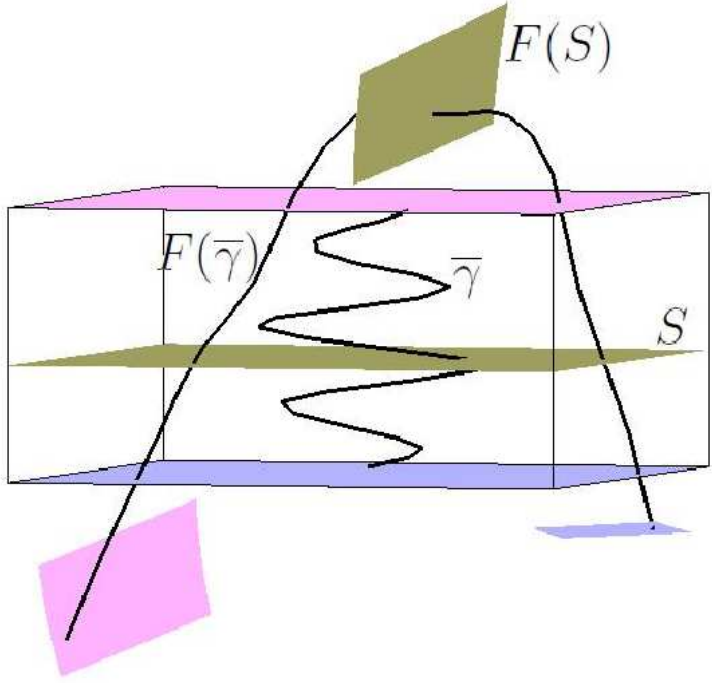}
\caption{With reference to the parallelepiped $\mathcal R$ in Figure 2,
reproduced here at a different scale, we show that the $F$-image set of
an arbitrary path $\gamma$ joining in $\mathcal R$ the two components of
the boundary set $\mathcal R^-$ intersects $\mathcal R$ twice. In particular, this is due
to the fact that the horizontal faces $\mathcal R^-_{\ell}$ and $\mathcal
R^-_{r}$ are mapped by $F$ below $\mathcal R^-_{\ell},$ in conformity with conditions $(C1)$ and $(C2),$  
and that the flat surface $S$ of the middle points w.r.t. the 
$z$-coordinate in $\mathcal R$ is mapped by $F$ above
$\mathcal R^-_{r},$ in agreement with condition $(C3^{'}).$ }
\end{figure}

\begin{figure}[ht]\label{fig-4}
\centering
\includegraphics[scale=0.7]{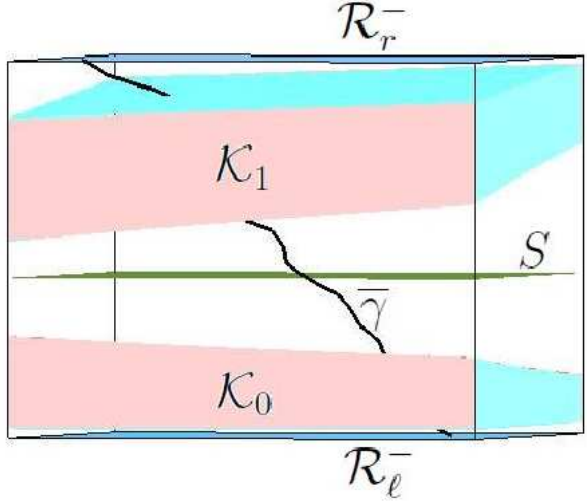}
\caption{Since $F(S)\cap\mathcal R=\emptyset$ (see Figure 4), then ${\mathcal R\cap F(\mathcal R)={\mathcal K}_0\cup{\mathcal K}_1},$ with 
${\mathcal K}_0$ and ${\mathcal K}_1$ disjoint.}
\end{figure}

\begin{figure}[ht]\label{fig-5}
\centering
\includegraphics[scale=1.2]{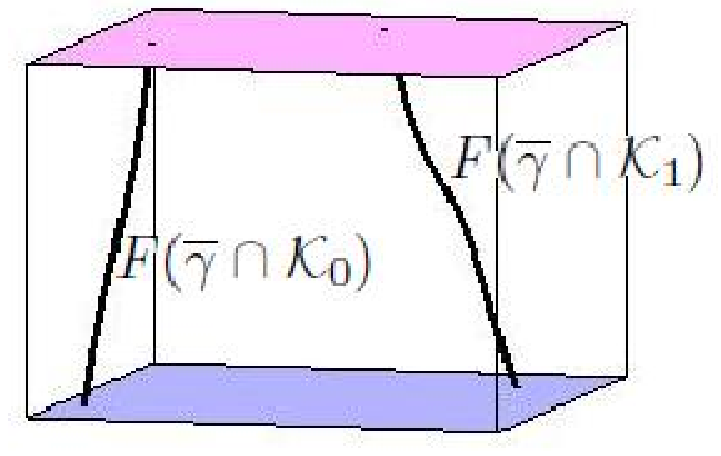}
\caption{Given the arbitrary path $\gamma$ in Figure $5$ joining in $\mathcal R$ the two components of
$\mathcal R^-,$ we show that the $F$-image sets of $\overline\gamma\cap{\mathcal K}_0$ and of $\overline\gamma\cap{\mathcal K}_1$ join $\mathcal R^-_{\ell}$ with $\mathcal R^-_{r},$ as required by the SAP property.}
\end{figure}

Before proving Theorem \ref{th1}, we make some comments on the conditions in $(H1)$--$(H5).$
First of all, notice that those conditions imply that $x_{\ell}+z_{\ell}>0$ and $x_{\ell}+y_{\ell}+z_{\ell}>0$ and thus there are no issues with the definition of $F$ on $\mathcal R$\,\,\footnote{Notice that, with our conditions on the parameters, it is immediate to check that also the functions we will introduce in the proof of Theorem \ref{th1} will be well defined, even when not explicitly remarked.}.
We also remark that we chose to split $(H1)$--$(H5)$ according to the corresponding conditions $(C1)$--$(C5)$ in the next proof they allow to verify.
Moreover we stress that the assumptions in $(H1)$--$(H5)$ are consistent, i.e., there exist parameter configurations satisfying them all. For instance, we checked that they are fulfilled for $c_1=0.4,\,c_2=0.55,\,c_3=0.6,\,\alpha=17,$ $x_{\ell}=0.5766666668,\,x_r=0.6316666668,\,y_{\ell}=0.3366666668,\,y_r=.04516666668,\,z_{\ell}=0,$ $z_r=0.3951779684.$ These are the same parameter values we used to draw Figures $2$--$6$, with the only exception of $z_{\ell}$ that in those pictures is slightly negative. Although this makes no sense from an economic viewpoint, as the variables $x,\,y$ and $z$ represent the output of the three firms, we made such choice in order to make the pictures easier to read. In fact, choosing $z_{\ell}=0,$ then 
$F(\mathcal R^-_{\ell})\subseteq \mathcal R^-_{\ell}$ and thus the crucial set $F(\mathcal R^-_{\ell})$ would have been not visible in Figures 2--4.
With this respect, we also remark that in Figure 3 the $x$-axis has been reversed in order to make the double folding of $F(\mathcal R)$ more evident.\\
In regard to the choice of the parameter values in \eqref{par}, as mentioned above, they are the same as in \cite{NaTr->}, except for $\alpha,$ which is larger here. In fact, numerical exercises we performed show that when $\alpha$ increases it becomes easier to find a domain where to apply the SAP technique. On the other hand, it seems not possible to apply our method for a sensibly smaller value of $\alpha.$
The impossibility of reducing $\alpha$ much below $17$ comes from the fact that, as it is immediate to verify, when such parameter decreases it becomes more and more difficult to have all conditions in $(H2)$ and $(H3)$ fulfilled and with $\alpha=10$ it seems just impossible. The situation would slightly improve dealing with 
$(C2)$ and $(C3)$ below, instead of $(C2)$ and $(C3^{'})$ as we actually do in order to simplify our argument, but still computer plots suggest it is not possible to have both conditions satisfied when $\alpha=8,$ that is the largest value considered in \cite{NaTr->}.

\medskip

\begin{proof}
We show that, for the parameter values in \eqref{par}, any choice of $x_i,\,y_i,\,z_i,\,i\in\{\ell,r\},$ fulfilling 
$(H1)$--$(H5)$ guarantees that the image under the map $F$ of any path
$\gamma=(\gamma_1,\gamma_2,\gamma_2):[0,1]\to\mathcal R=\mathcal R(x_i,y_i,z_i)$ joining the sets
$\mathcal R^-_{\ell}$ and $\mathcal
R^-_{r}$ defined in \eqref{eq-or} satisfies the following
conditions:
\begin{itemize}
\item [$(C1)$] $F_3(\gamma(0))\le z_{\ell}\,;$
\item [$(C2)$] $F_3(\gamma(1))\le z_{\ell}\,;$
\item [$(C3)$] $\exists \,\, t^*\in \, (0,1): F_3(\gamma(t^*))>z_{r}\,;$
\item [$(C4)$] $F_1(\gamma(t))\subseteq [x_{\ell},x_{r}],\,\forall t\in [0,1]\,;$
\item [$(C5)$] $F_2(\gamma(t))\subseteq [y_{\ell},y_{r}],\,\forall t\in [0,1]\,.$
\end{itemize}
Broadly speaking, conditions $(C1)$--$(C3)$ describe an expansion
with folding along the $z$--coordinate. In fact, the image $F\circ
\gamma$ of any path $\gamma$ joining in ${\mathcal R}$ the sides ${\mathcal R}^-_{\ell}$ and
${\mathcal R}^-_r$ crosses a first time the parallelepiped ${\mathcal
R}$ for $t\in(0,t^*)$ and then crosses ${\mathcal R}$ back
again for $t\in(t^*,1).$ Conditions $(C4)$ and $(C5)$ imply instead a
contraction along the $x$--coordinate and the $y$--coordinate, respectively.

\smallskip
\noindent
Actually, in order to simplify the exposition, instead of $(C3),$ we will check that the stronger condition
\begin{itemize}
\item [$(C3^{'})$] $F_3\left(x,y,\frac{z_{\ell}+z_{r}}{2}\right)>z_r,\,\forall (x,y)\in [x_{\ell},x_{r}]\times[y_{\ell},y_{r}],$
\end{itemize}
is satisfied, which means that the inequality in $(C3)$ holds for any
$t^*\in (0,1)$ such that
$\gamma(t^*)=\big(x,y,\frac{x_{\ell}+x_{r}}{2}\big),$ for some
$(x,y)\in [x_{\ell},x_{r}]\times[y_{\ell},y_{r}].$ Notice that 
\begin{equation}\label{s}
S:=\Big\{\Big(x,y,\frac{z_{\ell}+z_{r}}{2}\Big):(x,y)\in [x_{\ell},x_{r}]\times[y_{\ell},y_{r}]\Big\}\subseteq\mathcal R
\end{equation}
is the flat surface of middle points w.r.t. the 
$z$--coordinate in $\mathcal R$ depicted in Figure 4.\\
Setting
$$\mathcal R_0:=\Big\{(x,y,z)\in\mathbb R^3: (x,y)\in [x_{\ell},x_{r}]\times[y_{\ell},y_{r}],\,z\in\Big[z_{\ell},\frac{z_{\ell}+z_{r}}{2}\Big]\Big\},$$

$$\mathcal R_1:=\Big\{(x,y,z)\in\mathbb R^3: (x,y)\in [x_{\ell},x_{r}]\times[y_{\ell},y_{r}],\,z\in\Big[\frac{z_{\ell}+z_{r}}{2},z_r\Big]\Big\},$$

and
$$\mathcal K_0:=\mathcal R_0\cap F(\mathcal R) \quad \mbox{ and } \quad  \mathcal K_1:=\mathcal R_1
\cap F(\mathcal R)$$
(see Figure 5), we claim that $(C1),(C2),(C3^{'}),(C4)$ and
$(C5)$ together imply \eqref{eq-ks}. Notice at first that $\mathcal K_0$
and $\mathcal K_1$ are disjoint because, thanks to
condition $(C3^{'}),$ the set $S$ in \eqref{s} is mapped by $F$ outside $\mathcal R$ (see Figure 4). 
Furthermore, by $(C1),\,(C2)$ and
$(C3^{'}),$ for every path $\gamma: [0,1]\to {\mathcal R}$ such
that $\gamma(0)$ and $\gamma(1)$ belong to different components of ${\mathcal R}^-,$
there exist two disjoint
sub-intervals $[t_0',t_0''],\, [t_1',t_1'']\subseteq [0,1]$ such
that, setting
$\sigma_0:=\gamma|_{[t_0',t_0'']}:[t_0',t_0'']\to {\mathcal K}_0$ and 
$\sigma_1:=\gamma|_{[t_1',t_1'']}:[t_1',t_1'']\to {\mathcal K}_1,$
it holds that $F(\sigma_0(t_0'))$ and
$F(\sigma_0(t_0''))$ belong to different components of ${\mathcal R}^-,$ as well as $F(\sigma_1(t_1'))$ and
$F(\sigma_1(t_1'')).$ Moreover, from $(C4)$ and $(C5)$ it follows that 
$F(\sigma_0(t))\in {\mathcal R},\,\forall\, t\in [t_0',t_0'']$ and 
$F(\sigma_1(t))\in {\mathcal R},\,\forall\, t\in [t_1',t_1''].$\\
This means that $({\mathcal K}_i,F): {\widetilde{\mathcal R}} \stretchx {\widetilde{\mathcal R}},\,i=0,1,$ and
our claim is thus proved.\\
Once that the stretching condition in \eqref{eq-ks} is achieved,
the conclusion of the theorem follows by Theorem \ref{th-ch} \footnote{Notice that, by the choice of $\mathcal K_0$ and $\mathcal K_1,$ the invariant chaotic set $\mathcal I\subseteq\mathcal K_0\cup\mathcal K_1$
in Definition \ref{def-ch}
lies entirely in the first quadrant and therefore makes economic sense for the application in question.}.

\medskip

In order to complete the proof, let us verify that any choice of the parameters as in \eqref{par} and of the domain 
${\mathcal R}=\mathcal R(x_i,y_i,z_i)$ in agreement with $(H1)$--$(H5)$
implies that conditions $(C1),(C2),(C3^{'}),(C4)$ and
$(C5)$ are fulfilled for any path
$\gamma:[0,1]\to\mathcal R$ joining 
$\mathcal R^-_{\ell}$ and 
$\mathcal R^-_{r}$ \footnote{Just to fix the ideas, in what follows we will assume that $\gamma(0)\in\mathcal R^-_{\ell}$ and $\gamma(1)\in\mathcal R^-_{r}.$}.
In so doing, we will prove that the inequality in $(C1)$ is indeed an equality.\\
Let us start with the verification of $(C1).$ Since $F_3(x,y,z)=z\left(1-\alpha c_3+\frac{\alpha (x+y)}{(x+y+z)^2}\right)$
and $\gamma(0)\in\mathcal R^-_{\ell}= [x_{\ell},x_r]\times [y_{\ell},y_r]\times\{z_{\ell}\}=[x_{\ell},x_r]\times [y_{\ell},y_r]\times\{0\}$ by $(H1),$
it then follows that $\gamma_3(0)=0$ and thus
$0=F_3(\gamma(0))\le z_{\ell}=0,$ as desired.\\
In regard to $(C2),$ we have to verify that $F_3|_{\mathcal R^-_{r}}\le 0,$ that is, $F_3(x,y,z_r)\le 0,$ $\forall (x,y)\in [x_{\ell},x_{r}]\times[y_{\ell},y_{r}].$ 
Setting $A:=x+y,$ we consider, instead of $F_3|_{\mathcal R^-_{r}},$ the one-dimensional function\footnote{In several steps of the proof, 
instead of studying the original problem, through a substitution we will be lead to consider a lower dimensional one. Alternatively, we could use the Kuhn-Tucker Theorem for constrained maximization problems. We decided to follow the former approach because it is more elementary and requires less computations. However, we stress that the two approaches require to impose the same conditions $(H1)$--$(H5)$ on the parameters.}
$$\phi:[x_{\ell}+y_{\ell},x_r+y_r]\to\mathbb R,\quad \phi(A):=z_r\left(1-\alpha c_3+\frac{\alpha A}{(A+z_r)^2}\right).$$
Computing the first derivative of $\phi,$ we get
$\phi^{\,'}(A)=z_r\,\alpha\left(\frac{-A+z_r}{(A+z_r)^3}\right),$ which vanishes at $A=z_r.$ However, since by $(H2)$ we have $x_{\ell}+y_{\ell}>z_r,$ then $\phi{\,'}(A)<0,$ $\forall A\in[x_{\ell}+y_{\ell},x_r+y_r].$ Hence, $F_3|_{\mathcal R^-_{r}}\le F_3(x_{\ell},y_{\ell},z_r)$ and thus, in order to have $(C2)$ satisfied, it suffices that $F_3(x_{\ell},y_{\ell},z_r)\le 0.$ Imposing such condition, we find 
$z_r\left(1-\alpha c_3+\frac{\alpha (x_{\ell}+y_{\ell})}{(x_{\ell}+y_{\ell}+z_r)^2}\right)\le 0,$ which is fulfilled when
$\frac{\alpha c_3-1}{\alpha}\ge\frac{x_{\ell}+y_{\ell}}{(x_{\ell}+y_{\ell}+z_r)^2}.$ Making $z_r$ explicit, this 
holds when $z_r\ge\sqrt{\frac{\alpha}{\alpha c_3-1}(x_{\ell}+y_{\ell})}-(x_{\ell}+y_{\ell}),$ that is, when
$(H2)$ is fulfilled. Notice that the latter is a ``true'' restriction, since, still by $(H2),$ the right hand side of the above inequality is positive.
The verification of $(C2)$ is complete.\\
As regards $(C3^{'}),$ we need to check that $F_3\left(x,y,\frac{z_{\ell}+z_{r}}{2}\right)>z_r,\,\forall (x,y)\in [x_{\ell},x_{r}]\times[y_{\ell},y_{r}],$ that is, recalling the definition of $S$ in \eqref{s}, $F_3|_S>z_r.$ Notice that, by $(H1),$
$\frac{z_{\ell}+z_{r}}{2}=\frac{z_{r}}{2}.$
Analogously to what done above, instead of $F_3|_S,$ let us consider the one-dimensional function
$$\varphi:[x_{\ell}+y_{\ell},x_r+y_r]\to\mathbb R,\quad \varphi(A):=\frac{z_{r}}{2}\left(1-\alpha c_3+\frac{\alpha A}{\big(A+\frac{z_{r}}{2}\big)^2}\right).$$
Since $x_{\ell}+y_{\ell}>z_r>\frac{z_r}{2},$ by the previous analysis we know that $\varphi(A)\ge \varphi(x_r+y_r)=F_3\left(x_r,y_r,\frac{z_{r}}{2}\right).$
Hence, in order to have $F_3|_S>z_r,$ it suffices that $F_3\left(x_r,y_r,\frac{z_{r}}{2}\right)>z_r,$ that is, 
$$\frac{z_{r}}{2}\left(1-\alpha c_3+\frac{\alpha (x_r+y_r)}{\big(x_r+y_r+\frac{z_{r}}{2}\big)^2}\right)>z_r.$$ 
Since $z_{r}>0,$ making $z_r$ explicit, we find 
$$z_r<2\left(\sqrt{\frac{\alpha}{\alpha c_3+1}(x_{r}+y_{r})}-(x_{r}+y_{r})\right)$$
and this condition is satisfied thanks to $(H3).$ Hence $(C3)$ is verified.\\
In order to check $(C4),$ we need to show the two inequalities $F_1(x,y,z)\le x_r,$ $\forall (x,y,z)\in\mathcal R$ and 
$F_1(x,y,z)\ge x_{\ell},\,\forall (x,y,z)\in\mathcal R,$ which are satisfied if 
$${\displaystyle{\max_{(x,y,z)\in\mathcal R}F_1(x,y,z)}}\le x_r \qquad \mbox{and} \qquad {\displaystyle{\min_{(x,y,z)\in\mathcal R}F_1(x,y,z)}}\ge x_{\ell}\,,$$ 
respectively\footnote{Notice that such maximum and minimum values exist by the Weierstrass Theorem.}.\\
Instead of considering $F_1|_{\mathcal R},$ setting $B:=y+z$ and $T:=[x_{\ell},x_r]\times[y_{\ell}+z_{\ell},y_r+z_r],$ we deal with the bidimensional function
$$\Phi:T\to\mathbb R,\quad \Phi(x,B):=\frac{2x+B-c_1(x+B)^2}{2}\,,$$
whose partial derivatives are
$$\frac{\partial\Phi}{\partial x}=1-c_1(x+B)\qquad \mbox{ and } \qquad \frac{\partial\Phi}{\partial B}=\frac{1}{2}-c_1(x+B).$$
Since they do not vanish contemporaneously, there are no critical points in the interior of $T.$ We then study 
$\Phi$ on the boundary of its domain.\\ 
As concerns $\Phi_1(B):=\Phi|_{\{x_{\ell}\}\times[y_{\ell}+z_{\ell},y_r+z_r]}(x,B)=\Phi(x_{\ell},B),$ we have that 
$\Phi_1^{\,'}(B)=\frac{1}{2}-c_1(x_{\ell}+B),$ which vanishes at $\overline B=\frac{1}{2c_1}-x_{\ell}.$ This is the maximum point of $\Phi_1$ if $\overline B\in [y_{\ell}+z_{\ell},y_r+z_r].$ But that is guaranteed by the conditions in $(H4).$\\
Similarly, setting $\Phi_2(B):=\Phi|_{\{x_{r}\}\times[y_{\ell}+z_{\ell},y_r+z_r]}(x,B)=\Phi(x_{r},B),$ we find that its maximum point, still by $(H4),$ is given by $\widehat B=\frac{1}{2c_1}-x_{r}\in [y_{\ell}+z_{\ell},y_r+z_r].$\\
In regard to $\Phi_3(x):=\Phi|_{[x_{\ell},x_r]\times\{y_{\ell}+z_{\ell}\}}(x,B)=\Phi(x,y_{\ell}+z_{\ell}),$ we have $\Phi_3^{\,'}(x)=1-c_1(x+y_{\ell}+z_{\ell}),$  
which vanishes at $\overline x=\frac{1}{c_1}-(y_{\ell}+z_{\ell}).$ By the conditions in $(H4),$ $\overline x>x_r$ and thus $\Phi_3(x)$ is increasing on $[x_{\ell},x_r]$. Analogously, since $\widehat x=\frac{1}{c_1}-(y_r+z_r)>x_r,$ it holds that $\Phi_4(x):=\Phi|_{[x_{\ell},x_r]\times\{y_r+z_r\}}(x,B)=\Phi(x,y_r+z_r)$ is increasing on $[x_{\ell},x_r].$ Summarizing, the two candidates for the maximum point of $\Phi$ on $T$ are $\big(x_{\ell},\frac{1}{2c_1}-x_{\ell}\big)$ and $\big(x_{r},\frac{1}{2c_1}-x_{r}\big).$ A direct computation shows that $\Phi\big(x_{\ell},\frac{1}{2c_1}-x_{\ell}\big)<\Phi\big(x_{r},\frac{1}{2c_1}-x_{r}\big),$ and thus 
${\displaystyle{\max_{(x,y,z)\in\mathcal R}F_1(x,y,z)}}=\Phi\big(x_{r},\frac{1}{2c_1}-x_{r}\big).$ Hence, it is now easy to verify that the inequality ${\displaystyle{\max_{(x,y,z)\in\mathcal R}F_1(x,y,z)}}\le x_r$ is satisfied when $x_r\ge \frac{1}{4c_1},$ the latter being among the assumptions in $(H4).$\\
The analysis above also suggests that the two candidates for the minimum point of $\Phi$ on $T$ are $(x_{\ell},y_{\ell}+z_{\ell})$ and $(x_{\ell},y_r+z_r).$ Straightforward calculations show that, if $x_{\ell}\le \frac{1}{2c_1}\left(1-c_1(y_{\ell}+y_{r}+z_{\ell}+z_{r})\right),$ then 
$\Phi(x_{\ell},y_{\ell}+z_{\ell})\le \Phi(x_{\ell},y_r+z_r).$ Hence, again by $(H4),$ 
${\displaystyle{\min_{(x,y,z)\in\mathcal R}F_1(x,y,z)}}=\Phi\big(x_{\ell},y_{\ell}+z_{\ell}\big).$ The inequality ${\displaystyle{\min_{(x,y,z)\in\mathcal R}F_1(x,y,z)}}\ge x_{\ell}$ is thus satisfied when $\sqrt{\frac{y_{\ell}+z_{\ell}}{c_1}}-(y_{\ell}+z_{\ell})\ge x_{\ell},$ which is among the conditions in $(H4).$\\
This concludes the verification of $(C4).$\\
Let us finally turn to $(C5).$ In order to check it, we have to show that 
\begin{equation}\label{mm}
{\displaystyle{\max_{(x,y,z)\in\mathcal R}F_2(x,y,z)}}\le y_r \qquad \mbox{and} \qquad {\displaystyle{\min_{(x,y,z)\in\mathcal R}F_2(x,y,z)}}\ge y_{\ell}\,.
\end{equation} 
Instead of $F_2|_{\mathcal R},$ setting $D:=x+z,$ we deal with the one-dimensional function
$$\psi:[x_{\ell}+z_{\ell},x_r+z_{r}]\to\mathbb R,\quad \psi(D):=\sqrt{\frac{D}{c_2}}-D,$$
whose derivative is $\psi{'}(D)=\frac{1}{2\sqrt{c_2 D}}-1.$ It vanishes at $\overline D=\frac{1}{4c_2},$ which by $(H5)$ is smaller than $x_{\ell}+z_{\ell}.$ Thus 
${\displaystyle{\max_{(x,y,z)\in\mathcal R}F_2(x,y,z)}}=\psi(x_{\ell}+z_{\ell})$ and ${\displaystyle{\min_{(x,y,z)\in\mathcal R}F_2(x,y,z)}}=\psi(x_{r}+z_{r}).$ Hence, the first condition in \eqref{mm} is satisfied if $\psi(x_{\ell}+z_{\ell})\le y_r$ and the second condition  
is fulfilled if $\psi(x_{r}+z_{r})\ge y_{\ell}.$ It is easy to see that both inequalities are fulfilled thanks to $(H5)$ and this concludes the verification of $(C5).$\\
The proof is complete.
\end{proof}

\section{Conclusions}\label{sec-a3}
$~$

In this paper we have recalled what the SAP method consists in and we have applied that topological technique to rigorously prove the existence of chaotic sets for the triopoly game model in \cite{NaTr->}. By ``chaotic sets'' we mean invariant domains on which the map describing the system under consideration is semiconjugate to the Bernoulli shift (implying the features in Remark \ref{cons}) and where periodic points are dense.
However, we stress that we did not say anything about the attractivity of those chaotic sets. In fact, in general, the SAP method does not allow to draw any conclusion in such direction.
For instance, when performing numeric simulations for the parameter values in \eqref{par},
no attractor appears on the computer screen. The same issue emerged with the bidimensional models
considered in \cite{MePiZa-09}. The fact that the chaotic set is repulsive can be a good signal as regards the overlapping generations model therein, for 
which we studied a backward moving system, since the forward moving one was defined only implicitly and it was not possible to
invert it. Indeed, as argued in \cite{MePiZa-09}, a repulsive chaotic set for the backward moving system possibly gets transformed 
into an attractive one for a related forward moving system through Inverse Limit Theory (ILT).
In general, however, one just deals with a forward moving dynamical system and this kind of argument cannot be employed. 
For instance, both in the duopoly game model in \cite{MePiZa-09} and in the triopoly game model analyzed in the present paper, we are able to
prove the presence of chaos for the same parameter values considered in the literature, except for a bit larger speed of adjustment $\alpha.$ It makes economic sense that complex dynamics arise when firms are more reactive, but unfortunately for such parameter values no chaotic attractors can be found
via numerical simulations.\\
What we want to stress is that this is not a limit of the SAP method: such issue is instead related to the possibility of performing computations by hands. To see what is the point, let us consider the well-known case of the logistic map $f:[0,1]\to\mathbb R,$ $f(x)=\mu x(1-x),$ with $\mu>0.$
As observed in \cite{MePiZa-09}, if we want to show the presence of chaos for it via the SAP method by looking at the first iterate, then we need $\mu>4.$ In this case, however, the interval $[0,1]$ is not mapped into itself and for almost all initial points in $[0,1]$ forward iterates limit to $-\infty.$
If we consider instead the second iterate, then the SAP method may be applied for values less than $4,$ for which chaotic attractors do exist. Figure 6 shows a possible choice for the compact sets $\mathcal K_0$ and $\mathcal K_1$ (denoted in the picture by $I_0$ and $I_1,$ since they are intervals) for the stretching relation to be satisfied when $\mu\sim3.88.$
\begin{figure}[ht]\label{fig-6}
\centering
\includegraphics[scale=0.25]{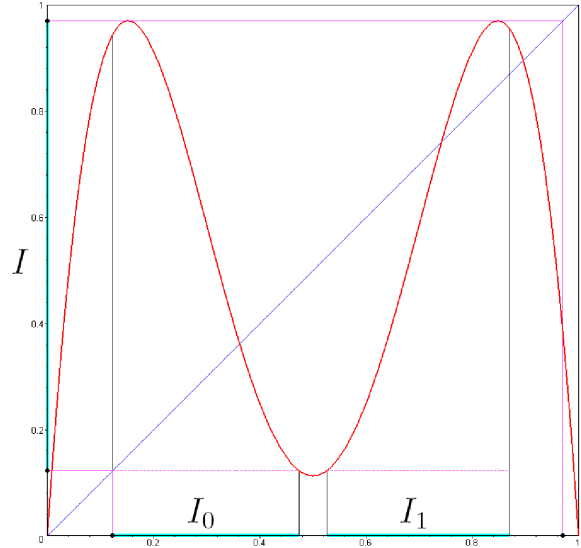}
\caption{The graph of the second iterate of the logistic map with $\mu\sim3.88.$}
\end{figure}
\noindent
This simple example aims to suggest that working with higher iterates may allow to reach an agreement between the conditions needed to employ the SAP method and those to find chaotic attractors via numerical simulations.\\
A possible direction of future study can then be the study of economically interesting but simple enough models, so that it is
possible to deal with higher iterates, in the attempt of rigorously proving the presence of chaos via the SAP technique for parameter values for which also computer simulations indicate the same kind of behavior.\\
Still in regard to chaotic attractors, we have observed that the SAP method works well for models presenting H\'enon-like attractors, due to the presence of a double folding, in turn related to the geometry required to apply our technique. On the other hand, a preliminary analysis seems to suggest that the SAP method is not easily applicable to models presenting a Neimark-Sacker bifurcation leading to chaos. A more detailed investigation of such kind of framework will be pursued, as well.\\
A further possible direction of future study is the analysis of continuous-time economic models with our technique, maybe in the context of   
LTMs, for systems switching between two different regimes, such as gross complements and gross substitutes.

\bigskip
\noindent
\textbf{Acknowledgements}.
Many thanks to Dr. Naimzada, Prof. Pini and Prof. Zanolin for useful discussions during the preparation of the paper.

\end{document}